\documentstyle[12pt]{article}
\if@twoside  m
\else
\fi
 \newcommand{\ie}{{\em i.e.}}
 \newcommand{\eg}{{\em e.g.}}
 
 \newcommand{\lhs}{{\em lhs }}
 \newcommand{\BB}{{\cal B}}
 \newcommand{\CC}{{\cal C}}
 \newcommand{\II}{{\cal I}}
 \newcommand{\KK}{{\cal K}}
 \newcommand{\LL}{{\cal L}}
 \newcommand{\NN}{{\cal N}}
 \newcommand{\VV}{{\cal V}}
 \newcommand{\XX}{{\cal X}}
 \newcommand{\R}{I\!\!R}
 \newcommand{\Z}{Z\!\!\!Z}
 \newcommand{\im}{{\rm Im\,}}
 \newcommand{\QED}{\mbox{\rule[-1.5pt]{6pt}{10pt}}}
 \newtheorem{claim}{Claim}[section]
 \newtheorem{theorem}[claim]{Theorem}
 \newtheorem{example}[claim]{Example}

\begin{document}

\vspace*{5mm}

\noindent
{\Large\bf A duality between Schr\"odinger operators \\ on graphs and
certain Jacobi matrices}
\vspace{5mm}

\begin{quote}
{\large P.~Exner}
\vspace{10mm}

{\em Nuclear Physics Institute, Academy of Sciences,
25068 \v Re\v z near Prague, \\
and Doppler Institute, FNSPE, Czech Technical University,
B\v rehov{\'a} 7, \\ 11519 Prague, Czech Republic \\
\rm exner@ujf.cas.cz}
\vspace{8mm}

The known correspondence between the Kronig--Penney
model and certain Jacobi matrices is extended to a wide class of
Schr\"odinger operators on graphs. Examples include rectangular
lattices with and without a magnetic field, or comb--shaped graphs
leading to a Maryland--type model.
\end{quote}

\section{Introduction}

Schr\"odinger operators on $\,L^2(\Gamma)\,$, where $\,\Gamma\,$ is
a graph, were introduced into quantum mechanics long time ago \cite{RuS}.
In recent years we have wittnessed of a renewed interest to them ---
see \cite{ARZ,GP,ES,BT,Ad,AL,AEL,GLRT,E2,E3,EG} and other, often
nonrigorous, studies quoted in these papers --- motivated mostly by the
fact that they provide a natural, if idealized, model of semiconductor
``quantum wire" structures.

Jacobi matrices, on the other hand, attracted a lot of attention in
the last decade, in particular, as a laboratory for random and almost
periodic systems. The most popular examples are the Harper and related
almost Mathieu equation dating back to \cite{Ha,Az,Ho}; for more recent
results and an extensive bibliography see, \eg,
\cite{Si,CFKS,AGHH,Be,La,Sh}. The underlying lattices are mostly
periodic of dimension one or two; however, more complicated examples
have also been studied \cite{Ma,Su}.

In case when $\,\Gamma\,$ is a line with an array of point interactions,
\ie, the Schr\"odinger operator in question is a Kronig--Penney--type
Hamiltonian, there is a bijective correspondence --- dubbed the {\em
French connection} by B.~Simon \cite{Si} --- between such systems and
certain Jacobi matrices \cite{AGHH,Ph,BFLT,DSS,GH,GHK}. The aim of
this paper is to show that the same duality can be established for a
wide class of Schr\"odinger operators on graphs, including the case
of a nonempty boundary. In general, the resulting Jacobi matrices
exhibit a varying ``mass".

\section{Preliminaries}

Let $\,\Gamma\,$ be a connected graph consisting of at most
countable families of {\em vertices} $\,\VV=\{\XX_j:\,j\in I\}\,$
and {\em links (edges)} $\,\LL=\{\LL_{jn}:\,
(j,n)\in I_{\LL}\subset I\times I\}\,$. We suppose that each pair of
vertices is connected by not more than one link. The set $\,\NN(\XX_j)=
\{\XX_n:\, n\in\nu(j)\subset I\setminus\{j\}\}\,$ of {\em neighbors}
of $\,\XX_j\,$, \ie, the vertices connected with $\,\XX_j\,$ by a single
link, is nonempty by hypothesis. The graph {\em boundary} $\,\BB\,$
consists of vertices having a single neighbor; it may be empty. We denote
by $\,I_\BB\,$ and $\,I_\II\,$ the index subsets in $\,I\,$
corresponding to $\,\BB\,$ and the graph {\em interior}
$\,\II:=\VV\setminus\BB\,$, respectively.

$\,\Gamma\,$ has a natural ordering by inclusions between vertex
subsets. We assume that it has a {\em local} metric structure in the
sense that each link $\,\LL_{jn}\,$ is isometric with a line segment
$\,[0,\ell_{jn}]\,$. The graph can be also equipped with a {\em global}
metric, for instance, if it is identified with a subset of $\,\R^\nu$.
Of course, the two metrics may differ at a single link.
Using the local metric, we are able to introduce the Hilbert space
$\,L^2(\Gamma):= \bigoplus_{(j,n)\in I_\LL} L^2(0,\ell_{jn})\,$; its
elements will be written as $\,\psi=\{\psi_{jn}:\, (j,n)\in I_\LL\}\,$
or simply as $\,\{\psi_{jn}\}\,$. Given a family of functions
$\,V:=\{V_{jn}\}\,$ with $\,V_{jn}\in L^{\infty}(0,\ell_{jn})\,$,
we define the operator $\,H_\CC\equiv H(\Gamma,\CC,V)\,$ by
   \begin{equation} \label{graph SO}
H_\CC\{\psi_{jn}\} \,:=\, \{\, -\psi''_{jn}+V_{jn}\psi_{jn}\,:\;
(j,n)\in I_\LL\,\}
   \end{equation}
with the domain consisting of all $\,\psi\,$ with $\,\psi_{jn}\in
W^{2,2}(0,\ell_{jn})\,$ subject to a set $\,\CC\,$ of boundary
conditions at the vertices which connect the boundary values
   \begin{equation} \label{boundary values}
\psi_{jn}(j):=\lim_{x\to 0+} \psi_{jn}(x)\,, \qquad
\psi'_{jn}(j):= \lim_{x\to 0+} \psi'_{jn}(x)\;;
   \end{equation}
we identify here the point $\,x=0\,$ with $\,\XX_j\,$.

There are many ways how to make the operator (\ref{graph SO})
self--adjoint; by standard results \cite{AGHH,RS} each vertex
$\,\XX_j\,$ may support an $\,N_j^2$ parameter family of boundary
conditions, where $\,N_j:= {\rm card\,} \nu(j)\,$. These self--adjoint
extensions were discussed in detail in \cite{ES}; here we restrict
ourselves to one of the following two possibilities which represent
in a sense extreme cases \cite{E2,E3}: \\
(a) {$\,\delta$ \em coupling:} at any $\,\XX_j\in\II\,$
we have $\;\psi_{jn}(j)=\psi_{jm}(j)=:\psi_j\,$ for all
$\,n,m\in\nu(j)\,$, and
   \begin{equation} \label{delta}
\sum_{n\in\nu(j)} \psi'_{jn}(j)\,=\, \alpha_j\psi_j
   \end{equation}
for some $\,\alpha_j\in\R\,$.  \\
(b) {$\,\delta'_s$ \em coupling:} $\;\psi'_{jn}(j)=
\psi'_{jm}(j)=:\psi'_j\,$ for all $\,n,m\in\nu(j)\,$, and
   \begin{equation} \label{delta'}
\sum_{n\in\nu(j)} \psi_{jn}(j)\,=\, \beta_j\psi'_j
   \end{equation}
for $\,\beta_j\in\R\,$. \\
The relations (\ref{delta}) and (\ref{delta'}) are independent of
$\,V\,$, since potentials are supposed to be essentially bounded.
At the graph boundary we employ the usual conditions,
   \begin{equation} \label{free end bc}
\psi_j \cos\omega_j+\psi'_j\sin\omega_j\,=\,0\,,
   \end{equation}
which can be written in either form (infinite values allowed);
for the sake of brevity we denote the set $\,\CC\,$ as
$\,\alpha=\{\alpha_j\}\,$ or $\,\beta=\{\beta_j\}\,$, respectively,
with the index running over $\,\VV\,$, and use the shorthands
$\,H_\alpha:= H(\Gamma,\{\alpha_j\},V)\,$ and  $\,H_\beta:=
H(\Gamma,\{\beta_j\},V)\,$.
Looking for solutions to the equations
   \begin{equation} \label{local SO}
H_\CC\psi\,=\, k^2\psi\,, \quad \CC=\alpha, \beta,
   \end{equation}
we shall consider the class $\,D_{loc}(H_\CC)\,$ which is
the subset in $\,\bigvee_{(j,n)\in I_\LL} L^2(0,\ell_{jn})\;$ (the
direct sum) consisting of the functions which satisfy all requirements
imposed at $\,\psi\in D(H_\CC)\,$ except the global square integrability.

The conditions (\ref{delta}) and (\ref{delta'}) define self--adjoint
operators also if the coupling constants are formally put equal to
infinity. We exclude this possibility, which corresponds, respectively,
to the Dirichlet and Neumann decoupling of the operator at $\,\XX_j\,$
turning the vertex effectively into $\,N_j\,$ points of the boundary.
On the other hand, we need it to state the result. Let
$\,H_\alpha^D\,$ and $\,H_\beta^N\,$ be the operators obtained from
$\,H_\alpha,\,H_\beta\,$ by changing the conditions (\ref{delta}),
(\ref{delta'}) at the points of $\,\II\,$ to Dirichlet
and Neumann, respectively, while at the boundary they are kept fixed.
We define $\,\KK_\alpha:=\{ k:\, k^2\in \sigma(H_\alpha^D),\:
\im k\ge 0\,\}\,$ and $\,\KK_\beta\,$ in a similar way.

\section{Main result}

Consider the operators $\,H_\alpha,\, H_\beta\,$ defined above.
We shall adopt the following assumptions:
   \begin{description}
   \vspace{-.8ex}
\item{\em (i)$\:$} There is $\,C>0\,$ such that $\,\|V_{jn}\|_{\infty}
\le C\,$ for all $\,(j,n)\in I_\LL\,$.
\vspace{-.8ex}
\item{\em (ii)} $\,\ell_0:=\inf \{\,\ell_{jn}:\, (j,n)\in I_\LL\}>0\,$.
\vspace{-.8ex}
\item{\em (iii)} $\,L_0:=\sup \{\,\ell_{jn}:\, (j,n)\in I_\LL\}<\infty\,$.
\vspace{-.8ex}
\item{\em (iv)} $\;N_0:=\max \{\,{\rm card\,} \nu(j):\, j\in I\,\}
<\infty\,$.
   \end{description}
On $\,\LL_{nj}\equiv [0,\ell_{jn}]\,$ (the right endpoint identified
with $\,\XX_j\,$) we shall denote as $\,u_{jn}^\CC,\; \CC=\alpha, \beta\,$,
the solutions to $\,-f''+V_{jn}f=k^2f\,$ which satisfy the boundary
conditions
$$
u_{jn}^\alpha(\ell_{jn})= 1\!-\!(u_{jn}^\alpha)'(\ell_{jn})=0\,,
\qquad 1\!-\!u_{jn}^\beta(\ell_{jn})= (u_{jn}^\beta)'(\ell_{jn})=0
$$
and
   \begin{eqnarray*}
v_{jn}^\alpha(0)= 1\!-\!(v_{jn}^\alpha)'(0)=0\,,
\quad 1\!-\!v_{jn}^\beta(0)= (v_{jn}^\beta)'(0)=0 \qquad & {\rm if} &
n\in I_\II\,,  \\
v_{jn}^\CC(0)=\sin\omega_n\,, \qquad
(v_{jn}^\CC)'(0)=-\cos\omega_n \qquad & {\rm if} &
n\in I_\BB\,;
   \end{eqnarray*}
their Wronskians are
$$
W_{jn}^\alpha= -v_{jn}^\alpha(\ell_{jn})\,,\qquad
W_{jn}^\beta= (v_{jn}^\beta)'(\ell_{jn})\,,
$$
respectively, or
   \begin{eqnarray*}
W_{jn}^\alpha= u_{jn}^\alpha(0)\,,
\quad W_{jn}^\beta= -(u_{jn}^\beta)'(0) \qquad & {\rm if} &
n\in I_\II\,,  \\
W_{jn}^\CC=-u_{jn}^\CC(0) \cos\omega_n
-(u_{jn}^\CC)'(0)\sin\omega_n \qquad & {\rm if} &
n\in I_\BB\,.
   \end{eqnarray*}
If not necessary we do not mark explicitly the dependence of these
quantities on $\,k\,$.

   \begin{theorem}
(a) Let $\,\psi\in D_{loc}(H_\alpha)\,$ solve (\ref{local SO})
for some $\,k\not\in\KK_\alpha\,$ with $\,k^2\in\R$, $\,\im k\ge 0\,$.
Then the corresponding boundary values (\ref{boundary values}) satisfy
the equation
   \begin{equation} \label{discrete delta}
\sum_{n\in\nu(j)\cap I_\II} {\psi_n\over W_{jn}^\alpha}\,-\,
\left(\, \sum_{n\in\nu(j)} {(v_{jn}^\alpha)'(\ell_{jn})
\over W_{jn}^\alpha} -\alpha_j\, \right)\psi_j\,=\, 0\,.
   \end{equation}
Conversely, any solution $\,\{\psi_j:\, j\in I_\II\}\,$ to
(\ref{discrete delta}) determines a solution of (\ref{local SO}) by
   \begin{eqnarray}
\psi_{jn}(x)= {\psi_n\over W_{jn}^\alpha}\,u_{jn}^\alpha(x)
-\,{\psi_j\over W_{jn}^\alpha}\,v_{jn}^\alpha(x) \qquad
& {\rm if} & n\in \nu(j)\cap I_\II\,, \label{reconstruction i} \\
\psi_{jn}(x)= -\,{\psi_j\over W_{jn}^\alpha}\,v_{jn}^\alpha(x) \qquad
& {\rm if} & n\in \nu(j)\cap I_\BB\,. \label{reconstruction b}
   \end{eqnarray}
(b) For a solution $\,\psi\in D_{loc}(H_\beta)\,$ of (\ref{local SO})
with $\,k\not\in\KK_\beta\,$, the above formulae are replaced by
   \begin{equation} \label{discrete delta'}
\sum_{n\in\nu(j)\cap I_\II} {\psi'_n\over W_{jn}^\beta}\,+\,
\left(\, \sum_{n\in\nu(j)} {v_{jn}^\beta(\ell_{jn})
\over W_{jn}^\beta} +\beta_j\, \right)\psi'_j\,=\, 0\,.
   \end{equation}
and
   \begin{eqnarray}
\psi_{jn}(x)= -\,{\psi'_n\over W_{jn}^\beta}\,u_{jn}^\beta(x)
+\,{\psi'_j\over W_{jn}^\beta}\,v_{jn}^\beta(x) \qquad
& {\rm if} & n\in \nu(j)\cap I_\II\,, \label{reconstruction' i} \\
\psi_{jn}(x)= \,{\psi'_j\over W_{jn}^\beta}\,v_{jn}^\beta(x) \qquad
& {\rm if} & n\in \nu(j)\cap I_\BB\,. \label{reconstruction' b}
   \end{eqnarray}
(c) Under (i), (ii), $\,\psi\in L^2(\Gamma)\,$ implies that the
solution $\,\{\psi_j\},\, \{\psi'_j\}\,$ of (\ref{discrete delta})
and (\ref{discrete delta'}), respectively, belongs to
$\,\ell^2(I_\II)\,$. \\
(d) The opposite implication is valid provided (iii), (iv) also
hold, and $\,k\,$ has a positive distance from from $\,\KK_\CC\,$.
   \end{theorem}
{\em Proof.} (a,b) We shall consider $\,H_\alpha\,$ throughout, the
argument for $\,H_\beta\,$ is analogous; for simplicity we drop the
superscript $\,\alpha\,$. If $\,n\in I_\II\,$, the transfer matrix on
$\,\LL_{nj}\,$ is
$$
T_{nj}(x,0)\,=\, W_{jn}^{-1} \left(
\begin{array}{cc} u_{jn}(x)-u'_{jn}(0)v_{jn}(x) & \; u_{jn}(0)v_{jn}(x)
\\ \\ u'_{jn}(x)-u'_{jn}(0)v'_{jn}(x) & \; u_{jn}(0)v'_{jn}(x) \end{array}
\right)\;;
$$
the Wronskian is nonzero for $\,k\not\in\KK_\alpha\,$. This yields an
expression of $\,\psi_{jn}(x)\,$ in terms of $\,\psi_{jn}(n)=:\psi_n\,$
and $\,\psi'_{jn}(n)\,$, in particular,
   \begin{eqnarray*}
\psi_j &\!\!:=\!\!& \psi_{jn}(j)= u'_{jn}(0)\psi_n+
v_{jn}(\ell_{jn})\psi'_{jn}(n)\,, \\ \\
-\psi'_{jn}(j) &\!=\!& {{1\!-\!u'_{jn}(0)v'_{jn}(\ell_{jn})}
\over W_{jn}} \psi_n+ v'_{jn}(\ell_{jn})\psi'_{jn}(n)\;;
   \end{eqnarray*}
the sign change at the \lhs of the last condition reflects the fact
that (\ref{boundary values}) defines the outward derivative at
$\,\XX_j\,$. We express $\,\psi'_{jn}(n)\,$ from the first relation
and substitute to the second one; this yields
$$
\psi'_{jn}(j)\,=\, -\,{\psi_n\over W_{jn}}\,+\, {v'_{jn}(\ell_{jn})
\over W_{jn}}\, \psi_j\,.
$$
If $\,n\in I_\BB\,$ we have instead
   \begin{eqnarray*}
&& \psi_j = u'_{jn}(0)\psi_n
-u_{jn}(0)\psi'_{jn}(n)\,, \\ \\
-\psi'_{jn}(j) &\!=\!& {{u'_{jn}(0)v'_{jn}(\ell_{jn})-\cos\omega_n}
\over W_{jn}}\, \psi_n +{{u_{jn}(0)v'_{jn}(\ell_{jn})\!-\!
\sin\omega_n}\over W_{jn}}\, \psi'_{jn}(n)\;;
   \end{eqnarray*}
we may write, of course, $\,\psi'_n:= \psi'_{jn}(n)\,$. Using the
relations $\,\psi_n\cos\omega_n +\psi'_n\sin\omega_n=0\,$, we find
$$
\psi'_{jn}(j)\,=\,{v'_{jn}(\ell_{jn}) \over W_{jn}}\, \psi_j\,,
$$
so (\ref{discrete delta}) follows from (\ref{boundary values}).
The transfer--matrix expression of $\,\psi_{jn}(x)\,$ together
with the mentioned formula for $\,\psi'_{jn}(n)\,$ yield
(\ref{reconstruction i}); (\ref{reconstruction b}) is checked in
the same way. \\
(c) Without loss of generality we may consider real $\,\psi\,$ only.
Any solution to the Schr\"odinger equation on $\,\LL_{nj}\,$ can
also be written as
$$
\psi_{jn}(x)\,=\, \psi_{jn}(n)w_{jn}(x) +\psi'_{jn}(n)v_{jn}(x)\,,
$$
where $\,w_{jn}\,$ is the normalized Neumann solution at $\,\XX_n\,$,
$\,w'_{jn}(0)= 1\!-\!w_{jn}(0)=0\,$. The argument of
\cite[Sec.III.2.1]{AGHH} cannot be used here, even in the classically
allowed region. Instead we employ a standard result of the
Sturm--Liouville theory \cite[Sec.1.2]{Mar} by which
   \begin{equation} \label{SL representation}
v_{jn}(x)\,=\, {\sin kx\over k}\,+\, \int_0^x\, K_D(x,y)\,
{\sin ky\over k}\, dy\,,
   \end{equation}
and similar representations are valid for $\,w_{jn}\,$ and their
derivatives; the kernels have an explicit bound in terms of
the potential $\,V_{jn}\,$. The latter is essentially bounded,
by (i) uniformly over $\,\LL\,$. Hence there is a positive
$\,\ell_1< {1\over 2}\,\ell_0\,$ such that
$$
\max\Big\lbrace\, |v_{jn}(x)|,\, |v'_{jn}(x)\!-\!1|,\,
|w_{jn}(x)\!-\!1|,\, |w'_{jn}(x)|\,\Big\rbrace\,<\, {1\over 10}
$$
for $\,x\in[0,\ell_1)\,$ and any $\,(j,n)\in I_\LL\,$; we infer that
   \begin{eqnarray*}
\psi_{jn}(x)^2 &\!=\!& \psi_{jn}(n)^2w_{jn}(x)^2\!
+\psi'_{jn}(n)^2v_{jn}(x)^2\! +2 \psi_{jn}(n)\psi'_{jn}(n)
w_{jn}(x)v_{jn}(x) \\ \\
&\!\ge\!& {9\over 20}\, \psi_{jn}(n)^2- {37\over 100}\,
\psi'_{jn}(n)^2\,.
   \end{eqnarray*}
An analogous estimate can be made for $\,\psi'_{jn}(n)^2\,$; summing
both of them we get
$$
\psi_{jn}(x)^2 +\psi'_{jn}(x)^2 \,\ge\, {2\over 25}\,
\left( \psi_{jn}(n)^2 +\psi'_{jn}(n)^2 \right)\,.
$$
If $\,\psi\in L^2(\Gamma)\,$, it belongs to $\,D(H_\alpha)\,$, so
$\,\psi'\in L^2(\Gamma)\,$ also holds. Let $\,\LL_{jn}(\ell_1)\,$
be the link $\,\LL_{jn}\,$ with the middle part
$\,(\ell_1,\ell_{jn}\!-\!\ell_1)\,$ deleted; then
   \begin{eqnarray*}
\|\psi\|_{L^2(\Gamma)}^2\!+ \|\psi'\|_{L^2(\Gamma)}^2
&\!=\!& \sum_{(j,n)\in I_\LL}
\int_{\LL_{jn}} (\psi_{jn}(x)^2\! +\psi'_{jn}(x)^2)\,dx \\ \\
&\!\ge\!& \sum_{(j,n)\in I_\LL} \int_{\LL_{jn}(\ell_1)}
(\psi_{jn}(x)^2\! +\psi'_{jn}(x)^2)\,dx\\ \\
&\!\ge\!& {4\ell_1\over 25}\, \sum_{(j,n)\in I_\LL}
(\psi_{jn}(n)^2\!+\psi'_{jn}(n)^2) \,\ge\, {4\ell_1\over 25}\,
\sum_{j\in I_\II} \psi_j^2\,,
   \end{eqnarray*}
which yields the result. \\
(d) First we need to show that $\,W_{jn}^{-1}$ has a uniform bound
for $\,k\not\in \overline\KK_\alpha\,$ Without loss of generality
we may suppose that $\,\inf\KK_\alpha>0\,$; otherwise we shift all
potentials by a constant. In view of the Sturm comparison theorem
\cite[Sec.V.6]{Re}, see also \cite[Secs.XIII.1,15]{RS}, the zeros
of $\,v_{jn}(\ell_{jn},\cdot)\,$ and $\,v'_{jn}(\ell_{jn},\cdot)\,$
form a switching sequence; hence
it is sufficient to prove that $\,{dW_{jn}(k)\over dk}\,$ is
bounded uniformly away from zero around each root $\,k_0\in
\KK_\alpha\,$ of $\,W_{jn}(k)=0\,$. Using the known identity
\cite[Sec.I.5]{Be}, \cite[Sec.XIII.3]{RS},
$$
{dv_{jn}(x,k)\over dk^2}\, v'_{jn}(x,k)\,-\,
{dv'_{jn}(x,k)\over dk^2}\, v_{jn}(x,k)\,=\,
\int_0^x\, v_{jn}(x,k)^2\, dx\,,
$$
at $\,x=\ell_{nj},\: k=k_0\,$, we find
$$
{dv_{jn}(\ell_{jn},k_0)\over dk}\,=\, {2k_0\over
v'_{jn}(\ell_{jn},k_0)}\:
\int_{\LL_{jn}} v_{jn}(x,k_0)^2\, dx\,.
$$
We have supposed that $\,k_0\ge\inf\KK_\alpha>0\,$; furthermore,
by the argument of the preceding part we have $\,\int_{\LL_{jn}}
v_{jn}(x,k_0)^2 dx \ge {81\over 200}\,\ell_1^2\,$. At the same
time, using the representation (\ref{SL representation}) we find
that $\,|v'_{jn}(\ell_{jn},k_0)|\le C_1\,$ for a positive $\,C_1\,$
independent of $\,j,n\,$. Hence to a given $\,k\not\in\overline
\KK_\alpha\,$ there is a $\,C_2>0\,$, again independent of $\,j,n\,$,
such that $\,W_{jn}\ge C_2\,$. The representation (\ref{SL
representation}) also yields another uniform bound,
$$
C_3\,:=\, \sup\Big\lbrace\, \|u_{jn}\|_{L^2},\, \|v_{jn}\|_{L^2},\,:
\; (j,n)\in I_\LL \,\Big\rbrace\,<\,\infty\,.
$$
The relations (\ref{reconstruction i}) and (\ref{reconstruction b})
together with (iv) now yield
   \begin{eqnarray*}
\|\psi\|^2 &\!=\!& \sum_{(j,n)\in I_\LL} \|\psi_{jn}\|^2 \,\le\,
\sum_{j,n\in I_\II} \, |W_{jn}|^{-2} \left\lbrace\,
|\psi_n|^2 \|u_{jn}\|^2\!+ |\psi_j|^2 \|v_{jn}\|^2 \right\rbrace \\ \\
&& +\, \sum_{n\in I_\BB} \, |\psi_j|^2\, {\|v_{jn}\|^2\over
|W_{jn}|^2}\,\le\,
2N_0 \left( C_3\over C_2 \right)^2 \sum_{j\in I_\II} |\psi_j|^2\,,
   \end{eqnarray*}
so the result follows. \quad \QED

\vspace{3mm}

\noindent
{\em Remarks.} 1. Since the operators $\,H_\CC\,$ are below bounded,
one can express for $\,k^2\in\rho(H_\CC)\cap \rho(H_0)\,$ the resolvent
difference $\,(H_\CC\!-k^2)^{-1}\!- (H_0\!-k^2)^{-1}$, where $\,H_0=
H^D_\alpha,\, H^N_\beta\,$, by Krein's formula \cite{Kr,Ne,BKN}; this
expression becomes singular under the condition (\ref{discrete delta})
or (\ref{discrete delta'}), respectively. \\
2. If two vertices are joined by more than a single link,
the theorem can be applied if a vertex with the free $\,\delta\,$
coupling, $\,\alpha_j=0\,$, is added at each extra link. Vertices
with just two neighbors are also useful, if we want to amend the
potentials $\,V_{jn}\,$ with a family of point interactions ---
for an example see \cite{E1}. \\
3. If (i)--(iv) are relaxed, the implications (b), (c) may still hold
with $\,\ell^2(I_\II)\,$ replaced by a weighted $\,\ell^2\,$ space.
If $\,\Gamma\,$ is equipped with a global metric, one can establish
a relation between the exponential decay of $\,\psi\,$ and of the
boundary--value sequences.

\section{Examples}

   \begin{example} (rectangular lattices): {\rm Let $\,\Gamma\,$
be a planar lattice graph whose basic cell is a rectangle of sides
$\,\ell_j,\: j=1,2\,$, and suppose that apart of the junctions,
the motion on $\,\Gamma\,$ is free, $\,V_{jn}=0\,$. By the proved
theorem, the equation (\ref{local SO}) leads in this case to
   \begin{eqnarray} \label{lattice Jacobi}
\lefteqn{(\phi_{n,m+1}+\phi_{n,m-1})\sin k\ell_1+
(\phi_{n+1,m}+\phi_{n-1,m})\sin k\ell_2} \nonumber \\ \\
&& \phantom{AAAAAA} + \left(V^\CC_{nm}(k)\mp
2\sin k(\ell_1\!+\ell_2)\right)
\phi_{nm}\,=\,0 \phantom{AAAAAA} \nonumber
   \end{eqnarray}
in the $\,\delta\,$ and $\,\delta'_s\,$ case, respectively, where
   \begin{equation} \label{rectangle potentials}
V^\alpha_{nm}(k)\,:=\,
-\,{\alpha_{nm}\over k}\,\sin k\ell_1 \sin k\ell_2\,, \qquad
V^\beta_{nm}(k)\,:=\,
-\,\beta_{nm}k\,\sin k\ell_1 \sin k\ell_2 \,,
   \end{equation}
and $\,\phi_{nm}= \psi_{nm},\, \psi'_{nm}\,$. In particular, for
a square lattice we get
$$
(h_0\psi)_{nm}- \left( 4\cos k\ell+ {\alpha_{nm}\over k}
\sin k\ell \right)\psi_{nm}\,=\,0
$$
and an analogous equation in the $\,\delta'_s\,$ case, where
$\,h_0\,$ is the conventional two--dimensional Laplacian
on $\,\ell^2(\Z^2)\,$. On the other hand, if $\,\ell_1\ne\ell_2\,$,
the free operator has a periodically modulated ``mass" which leads to
nontrivial spectral properties even if the underlying graph
Hamiltonian is periodic, \ie, $\,\alpha_{nm}=\alpha\,$ or
$\,\beta_{nm}=\beta\,$, with a non--zero coupling parameter. A
detailed discussion of this situation can be found in \cite{E2,E3,EG}.
}
   \end{example}

  \begin{example} (magnetic field added): {\rm Suppose that
$\,\Gamma\,$ is embedded into $\,\R^\nu$ in which there is a
magnetic field described by a vector potential $\,A\,$. The
boundary conditions (\ref{delta}) and (\ref{delta'}) are
modified replacing $\,\psi'_{jn}(j)\,$ by $\,\psi'_{jn}(j)
+iA_{jn}(j)\,$, where $\,A_{jn}(j)\,$ is the tangent component
of $\,A\,$ to $\,\LL_{jn}\,$ at $\,\XX_j\,$; we suppose
conventionally that $\,e=-1\,$. For the $\,\delta\,$
coupling this is well known \cite{ARZ}; in the $\,\delta'_s\,$
case one can check it easily.

There in no need to repeat the above argument, however, since
the magnetic case can be handled with the help of the unitary
operator $\,U:\: L^2(\Gamma)\to L^2(\Gamma)\,$ defined by
$$
(U\psi)_{jn}(x)\,:=\, \exp\left( i\int_{x_{jn}}^x A_{jn}(y)\,dy
\right)\, \psi_{jn}(x)\,,
$$
where $\,A_{jn}\,$ is again the tangent component of the vector
potential and $\,x_{jn}\,$ are fixed reference points. Then the
functions $\,(U\psi)_{jn}\,$ satisfy (\ref{delta}) and
(\ref{delta'}), respectively, and the equations (\ref{discrete
delta}) and (\ref{discrete delta'}) are replaced by
   \begin{eqnarray*}
\sum_{n\in\nu(j)\cap I_\II} {e^{iA_n}\over W_{jn}^\alpha}\,
\psi_n\,-\, \left(\, \sum_{n\in\nu(j)}
{(v_{jn}^\alpha)'(\ell_{jn}) \over W_{jn}^\alpha} -\alpha_j\,
\right)e^{iA_j} \psi_j &\!=\!& 0\,,
\\ \\
   \sum_{n\in\nu(j)\cap I_\II} {e^{iA_n}\over W_{jn}^\beta}\,
\psi'_n\,+\, \left(\, \sum_{n\in\nu(j)}
{v_{jn}^\beta(\ell_{jn}) \over W_{jn}^\beta} +\beta_j\,
\right) e^{iA_j}\psi'_j &\!=\!& 0\,,
   \end{eqnarray*}
provided the magnetic phase factors $\,A_j\,$ obey the consistency
conditions
$$
A_j\!-A_n\,=\, \int_{\LL_{jn}} A_{jn}(y)\, dy
$$
following from the continuity requirement.

In particular, if $\,\Gamma\,$ is a rectangular lattice considered
above, $\,A_{nm}:= \Phi(n\!-\!m)/2\,$ corresponds to a homogeneous
magnetic field in the circular gauge, with the flux $\,\Phi:=
B\ell_1\ell_2\,$ through a cell, and no other potential is present,
the equation (\ref{lattice Jacobi}) is replaced by
   \begin{eqnarray} \label{mg lattice Jacobi}
\lefteqn{\left(e^{i\Phi m/2}\phi_{n,m+1}
+e^{-i\Phi m/2}\phi_{n,m-1}\right)\sin k\ell_1+
\left(e^{-i\Phi n/2}\phi_{n+1,m}+
e^{i\Phi n/2}\phi_{n-1,m}\right)\sin k\ell_2} \nonumber \\ \\
&&
\phantom{AAAAAAA}+ \left(V^\CC_{nm}(k)\mp 2\sin k(\ell_1\!+\ell_2)\right)
\phi_{nm}\,=\,0\,, \phantom{AAAAAAAAAAAAAAA} \nonumber
   \end{eqnarray}
where the discrete potential is again given by (\ref{rectangle
potentials}). If $\,\ell_1=\ell_2\,$ and the coupling constants
vanish, $\,\alpha_{nm}= \beta_{nm}=0\,$, we obtain in this way the
discrete magnetic Laplacian of \cite{Sh} (or the Harper operator
of \cite{Ha,Be} for the Landau gauge), while for a general
rectangle the parameter $\,\lambda\,$ representing the
anisotropy in \cite{Sh} is replaced again by the periodically
modulated ``mass".
}
   \end{example}

  \begin{example} (comb--shaped graphs): {\rm Let $\,\Gamma\,$
consist of a line to which at the points $\,\XX_j:=jL\,$ line
segments of lengths $\,\ell_j\,$ are joined; at their ends we
impose the conditions (\ref{free end bc}). The equation
(\ref{local SO}) now yields
   \begin{equation} \label{comb Jacobi}
(h_0\phi)_j+ \left(V^\CC_j(k)\mp 2\cos kL\right)\phi_j\,=\,0
   \end{equation}
in the $\,\delta\,$ and $\,\delta'_s\,$ case, respectively, where
$\,\phi_j= \psi_j,\, \psi'_j\,$, and
$$
V^\alpha_j(k)\,:=\,-\left( {v'_j(\ell_j)\over v_j(\ell_j)}
+\alpha_j\right)\, {\sin kL\over k}\,, \quad
V^\beta_j(k)\,:=\,-\left( {v_j(\ell_j)\over v'_j(\ell_j)}
+\beta_j\right)\, k\,\sin kL\,.
$$
In particular, in the absence of an external potential on
$\,\Gamma\,$ the last relations become
$$
V^\alpha_j(k)=-\left( \cot(k\ell_j\!+\!\eta_j)
+{\alpha_j\over k}\right)\, \sin kL\,, \quad
V^\beta_j(k)=-\left( \tan(k\ell_j\!+\!\eta_j)
+\beta_j k\right)\, \sin kL\,,
$$
where $\,\eta_j:=\arctan(k\tan\omega_j)\,$. If the coupling
is ideal, $\,\alpha_j=0\,$ or $\,\beta_j=0\,$, the loose ends
correspond to the Dirichlet condition, $\,\omega_j\,=0\,$, and
the ``tooth" lengths are $\,\ell_j:=|j|\ell\,$, we get thus an
equation reminiscent of the  Maryland model \cite{CFKS,PGF} with
a fixed coupling strength and an additional periodic modulation
of the potential.
}
   \end{example}

\vspace{3mm}

\noindent
{\em Acknowledgments.} The author is grateful for the hospitality
extended to him in the Institute of Mathematics, University of Ruhr,
Bochum, where this work was done. The research has been partially
supported by the Grant AS No.148409 and the European Union Project
ERB--CiPA--3510--CT--920704/704.


\vspace{5mm}
\begin{thebibliography}{article}
   \bibitem{RuS}
K.~Ruedenberg, C.W.~Scherr: Free--electron network model for
conjugated systems, I.~Theory, {\em J.Chem.Phys.} {\bf 21} (1953),
1565--1581.
   \vspace{-1.8ex}
   \bibitem{ARZ}
J.E.~Avron, A.~Raveh, B.~Zur: Adiabatic transport in multiply
connected systems, {\em Rev.Mod.Phys.} {\bf 60} (1988), 873--915.
   \vspace{-1.8ex}
   \bibitem{GP}
N.I.~Gerasimenko, B.S.~Pavlov: Scattering problem on noncompact graphs,
{\em Teor.Mat.Fiz.} {\bf 74} (1988), 345--359.
   \vspace{-1.8ex}
   \bibitem{ES}
P.~Exner, P.~\v{S}eba: Free quantum motion on a branching graph, {\em
Rep. Math.Phys.} {\bf 28} (1989), 7-26.
   \vspace{-1.8ex}
   \bibitem{BT}
W.~Bulla, T.~Trenckler: The free Dirac operator on compact and non-compact
graphs, {\em J.Math.Phys.} {\bf 31} (1990), 1157-1163.
   \vspace{-1.8ex}
   \bibitem{Ad}
V.M.~Adamyan: Scattering matrices for microschemes, {\em Oper.Theory:
Adv. Appl.} {\bf 59} (1992), 1--10.
   \vspace{-1.8ex}
   \bibitem{AL}
Y.~Avishai, J.M.~Luck: Quantum percolation and ballistic conductance
on a lattice of wires, {\em Phys.Rev.} {\bf B45} (1992), 1074--1095.
   \vspace{-1.8ex}
   \bibitem{AEL}
J.E.~Avron, P.~Exner, Y.~Last: Periodic Schr\"odinger operators with
large gaps and Wannier--Stark ladders, {\em Phys.Rev.Lett.} {\bf 72}
(1994), 896--899.
   \vspace{-1.8ex}
   \bibitem{GLRT}
J.~Gratus, C.J.~Lambert, S.J.~Robinson, R.W.~Tucker: Quantum
mechanics on graphs, {\em J.Phys.} {\bf A27} (1994), 6881--6892.
   \vspace{-1.8ex}
   \bibitem{E2}
P.~Exner: Lattice Kronig--Penney models, {\em Phys.Rev.Lett.} {\bf 74}
(1995), 3503--3506.
   \vspace{-1.8ex}
   \bibitem{E3}
P.~Exner: Contact interactions on graphs, {\em submitted for publication}
   \vspace{-1.8ex}
   \bibitem{EG}
P.~Exner, R.~Gawlista: Band spectra of rectangular graph superlattices,
{\em in preparation}
   \vspace{-1.8ex}
   \bibitem{Ha}
P.G.~Harper: Single band motion of conduction electrons in a
uniform magnetic field, {\em Proc.Roy.Soc.(London)} {\bf A68} (1955),
874--878.
   \vspace{-1.8ex}
   \bibitem{Az}
M.Ya.~Azbel: Energy spectrum of a conduction electron in a magnetic field,
{\em Sov.Phys. JETP} {\bf 19} (1964), 634--647.
   \vspace{-1.8ex}
   \bibitem{Ho}
D.R.~Hofstadter:Energy levels and wavefunctions for Bloch electrons
in rational  and irrational magnetic fields,  {\em Phys.Rev.} {\bf B14}
(1976), 2239--2249.
   \vspace{-4.8ex}
   \bibitem{Si}
B.~Simon: Almost periodic Schr\"odinger operators: a review,
{\em Adv.Appl. Math.} {\bf 3} (1982), 463--490.
   \vspace{-1.8ex}
   \bibitem{CFKS}
H.L.~Cycon, R.G.~Froese, W.~Kirsch, B.~Simon: {\em Schr\"odinger
Operators}, Springer, Berlin 1987.
   \vspace{-1.8ex}
   \bibitem{AGHH}
S.~Albeverio, F.~Gesztesy, R.~H\o egh-Krohn, H.~Holden: {\em Solvable
Models in Quantum Mechanics}, Springer, Heidelberg 1988.
   \vspace{-1.8ex}
   \bibitem{Be}
J.~Bellissard: Gap labelling theorem for Schr\"odinger operators
in {\em Number Theory and Physics,} (M. ~Waldschmidt et al., eds.),
Springer, Heidelberg 1992; pp. 538--630.
   \vspace{-1.8ex}
   \bibitem{La}
Y.~Last: Zero measure spectrum for almost Mathieu operator, {\em
Commun. Math.Phys.} {\bf 164} (1994), 421--432.
   \vspace{-1.8ex}
   \bibitem{Sh}
M.A.~Shubin: Discrete magnetic Laplacian, {\em Commun.Math.Phys.}
{\bf 164} (1994), 259--275.
   \vspace{-1.8ex}
   \bibitem{Ma}
L.~Malozemov: The integrated density of states for the difference
Laplacian on the modified Koch curve, {\em Commun.Math.Phys.}
{\bf 156} (1993), 387--397.
   \vspace{-1.8ex}
   \bibitem{Su}
T.~Sunada: Generalized Harper operator on a graph, in {\em Workshop
on Zeta Function in Number Theory and Geometric Analysis in Honor
of Jun--Ichi Igusa}, Johns Hopkins University 1993.
   \vspace{-1.8ex}
  \bibitem{Ph}
P.~Phariseau: The energy spectrum of an amorhous substance, {\em Physica}
{\bf 26} (1960), 1185--1191.
   \vspace{-1.8ex}
   \bibitem{BFLT}
J.~Bellissard, A.~Formoso, R.~Lima, D.~Testard: Quasi--periodic
interactions with a metal--insulator transition, {\em Phys.Rev.}
{\bf B26} (1982), 3024--3030.
   \vspace{-1.8ex}
   \bibitem{DSS}
F.~Delyon, B.~Simon, B.~Souillard: From power pure point to
continuous spectrum in disordered systems, {\em Ann.Inst.
H.~Poincar\'e} {\bf A42} (1985), 283--309.
   \vspace{-1.8ex}
   \bibitem{GH}
F.~Gesztesy, H.~Holden: A new class of solvable models in quantum
mechanics describing point interactions on the line, {\em J.Phys.}
{\bf A20} (1987), 5157--5177.
   \vspace{-4.8ex}
   \bibitem{GHK}
F.~Gesztesy, H.~Holden, W.~Kirsch: On energy gaps in a new type of
analytically solvable model in quantum mechanics, {\em
J.Math.Anal.Appl.} {\bf 134} (1988), 9--29.
   \vspace{-1.8ex}
   \bibitem{RS}
M.~Reed, B.~Simon: {\em Methods of Modern Mathematical Physics,
II.~Fourier Analysis. Self--adjointness, IV. Analysis of Operators},
Academic Press, New York 1975, 1978.
   \vspace{-1.8ex}
   \bibitem{Mar}
V.A.~Marchenko: {\em Sturm--Liouville Operators and Applications},
Operator Theory: Advances and Applications, vol.22; Birkh\"auser,
Basel 1986.
   \vspace{-1.8ex}
   \bibitem{Re}
W.T.~Reid: {\em Ordinary Differential Equations}, J.Wiley, New York 1971.
   \vspace{-1.8ex}
   \bibitem{Kr}
M.G.~Krein: The theory of self--adjoint extensions of semibounded
Hermitian transformations and its applications~I,II,
{\em Mat.Sbornik}~{\bf 20} (1947), 431--495; {\bf 21} (1947), 365--404.
   \vspace{-1.8ex}
   \bibitem{Ne}
G.~Nenciu: To the theory of self--adjoint extensions of symmetric
operators with a spectral gap, {\em Funkc.Anal.Appl.}~{\bf 19} (1985),
81--82.
   \vspace{-1.8ex}
   \bibitem{BKN}
J.F.~Brasche, V.~Koshmanenko, H.~Neidhardt: New aspects of Krein's
extension theory, {\em Ukrainian J.Math.}~{\bf 46} (1994), 37--54.
   \vspace{-1.8ex}
  \bibitem{E1}
P.~Exner: The absence of the absolutely continuous spectrum for
$\,\delta'\,$ Wannier--Stark ladders, {\em J.Math.Phys.} {\bf 36} (1995),
to appear
   \vspace{-1.8ex}
   \bibitem{PGF}
R.~Prange, D.~Grempel, S.~Fishman: A solvable model of quantum
motion in an incommensurate potential,  {\em Phys.Rev.} {\bf B29}
(1984), 6500--6512.
   \vspace{-1.8ex}

\end{thebibliography}
\end{document}